\documentclass[doublespacing]{elsart}

\usepackage{graphics}
\usepackage{graphicx}
\usepackage{epsfig}

\usepackage{amssymb}
\usepackage{amsmath}
\usepackage{latexsym}

\begin{document}

\begin{frontmatter}


\title{The Slow Invariant Manifold\\ of the \\ Lorenz-Krishnamurthy Model}
%\author[Gi]{Jean-Marc Ginoux\corauthref{cor}},
%\corauth[cor]{Corresponding author.} \ead{ginoux@univ-tln.fr}
%\ead[url]{http://ginoux.univ-tln.fr}

%\address[Gi]{Laboratoire LSIS, CNRS, UMR 7297, Universit\'{e} de Toulon,\\ BP 20132, F-83957 La Garde cedex, France}

%\title{}

%\author[label1,label2]{}
%\address[label1]{}
%\address[label2]{}

%\author{}

%\address{}

\begin{abstract}

During this last decades, several attempts to construct \textit{slow invariant manifold} of the Lorenz-Krishnamurthy five-mode model of slow-fast interactions in the atmosphere have been made by various authors. Unfortunately, as in the case of many \textit{two-time scales singularly perturbed dynamical systems} the various asymptotic procedures involved for such a construction diverge. So, it seems that till now only the first-order and third-order approximations of this slow manifold have been analytically obtained. While using the \textit{Flow Curvature Method} we show in this work that one can provide the eighteenth-order approximation of the \textit{slow manifold} of the generalized Lorenz-Krishnamurthy model and the thirteenth-order approximation of the ``conservative'' Lorenz-Krishnamurthy model. The invariance of each \textit{slow manifold} is then established according to \textit{Darboux invariance theorem}.

\end{abstract}

\begin{keyword}
Lorenz-Krishnamurthy model; Flow Curvature Method; Darboux invariance; Fenichel theory; slow invariant manifold

%\PACS
\end{keyword}
\end{frontmatter}

% main text
\section{Introduction}
\label{intro}

The classical geometric theory developed originally by Andronov
\cite{1}, Tikhonov \cite{23} and Levinson \cite{15} stated that
\textit{singularly perturbed systems} possess \textit{invariant
manifolds} on which trajectories evolve slowly and toward which
nearby orbits contract exponentially in time (either forward and
backward) in the normal directions. These manifolds have been called
asymptotically stable (or unstable) \textit{slow manifolds}. Then,
Fenichel \cite{{5},{6},{7},{8}} theory\footnote{This theory was
independently established in Hirsch, \textit{et al.} \cite{12}} for
the \textit{persistence of normally hyperbolic invariant manifolds}
enabled to establish the \textit{local invariance} of \textit{slow
manifolds} that possess both expanding and contracting directions
and which were labeled \textit{slow invariant manifolds}.
During the last century, various methods have been developed in
order to determine the \textit{slow invariant manifold} analytical
equation associated with \textit{singularly perturbed systems}. The
seminal works of Wasow \cite{25}, Cole \cite{2}, O'Malley
\cite{{16},{17}} and Fenichel \cite{{5},{6},{7},{8}} to name but a
few, gave rise to the so-called \textit{Geometric Singular
Perturbation Method}. According to this theory, existence as well as
local invariance of the \textit{slow manifold} of \textit{singularly
perturbed systems} have been stated.\\ Then, the determination of the
\textit{slow manifold} analytical equation turned into a
regular perturbation problem \cite[p. 112]{16} in which one generally expected the
asymptotic validity of such expansion to breakdown \cite{17}.

\newpage

In recent publications a new approach to $n$--dimensional \textit{singularly perturbed systems} of ordinary differential
equations with two time scales, called \textit{Flow Curvature Method} \cite{{9},{9a},{9b},{9c},{9d},{9e}} has been developed. It consists in considering the \textit{trajectory curves} integral of such systems as curves in Euclidean $n$--space. Based on the use of local metrics properties of \textit{curvatures} inherent to \textit{Differential Geometry}, this method which does not require the use of asymptotic expansions, states that the
location of the points where the local \textit{curvature} of the \textit{trajectory curves} of such systems is null defines a $(n - 1)$--dimensional \textit{manifold} associated with this system and called \textit{flow curvature manifold}. The invariance of this manifold is then stated according to a theorem introduced by Gaston Darboux \cite{3} in 1878.

The laws governing the behavior of the atmosphere permit the simultaneous presence of oscillation modes such as \textit{quasi-geostrophic modes} and \textit{inertial-gravity modes}. The former which have periods of few days are generally referred as \textit{Rossby Waves} while the latter whose periods are of few hours are called \textit{Gravity Waves}. In 1980, a nine-dimensional primitive equation (PE) model of the atmosphere enabling to superpose \textit{Rossby} and \textit{Gravity Waves} was originally proposed by the late Edward Norton Lorenz \cite{15aa}. Few years later, Lorenz \cite{15a} simplified the nine-dimensional model to a five-dimensional model by truncating to just five modes: three \textit{Rossby Waves} coupled to two \textit{Gravity Waves}. This five-dimensional model can be considered as a \textit{two-time scales singularly perturbed dynamical system} with three \textit{slow variables} (\textit{Rossby Waves}) and two \textit{fast variables} (\textit{Gravity Waves}) in accordance with physical observations of atmospheric behavior. In numerical weather prediction a problem arises because raw fields data can not be used as initial conditions for such model, since even when the initial wind and pressure fields are both fairly realistic, \textit{Gravity Waves} will occur, if the fields are not in ``proper balance''. According to Camassa  \cite[p. 357]{2aa}: ``Small errors in the ``proper balance'' between these two time scales lead to abnormal evolution of gravity waves, which in turn causes appreciable deviation of weather forecasts from actual observation on the time scale of gravity waves.''. To solve this initialization problem, existence of a \textit{slow manifold}\footnote{This concept has been introduced by C. E. Leith \cite{15aaa} in 1980.}, consisting of trajectory curves (orbits) for which \textit{Gravity Waves} motion is absent, was first postulated for such model. Then, an iteration scheme was developed to find from the state (point in phase space) specified by field data a corresponding initial state on this \textit{slow manifold}, so that weather forecasts with these initial states can be accurate on the same time scale as \textit{Rossby Waves}. In their paper, Lorenz and Krishnamruthy \cite{15b} identified the variables representing \textit{Gravity Wave} activity as the ones which can exhibit \textit{fast} oscillations, and defined the \textit{slow manifold} as an invariant manifold in the five dimensional phase space for which \textit{fast} oscillations never develop. However, in a subsequent paper, Lorenz and Krishnamruthy \cite{15b} identified a trajectory curve (orbit) which by construction has to lie on the slow manifold. They followed its evolution numerically to show that sooner or later \textit{fast} oscillations developed, thereby implying, as pointed out by the title of their article, that a \textit{slow manifold} according to the definition does not exist for this model. Such a result gave rise to a series of articles published from 1991 to 1996 by S.J. Jacobs \cite{12a}, J.P. Boyd \cite{1a}, A.C. Fowler and G. Kember \cite{8a}, R. Camassa and Siu-Kei Tin \cite{2a} in which their authors proved the existence of a slow manifold for the Lorenz-Krishnamurthy model (LK). These latter proposed a \textit{generalized LK model}. More recently, M. Phani Sudheer, Ravi S. Nanjundiah  A.S. Vasudeva Murthy \cite{17a} ``revisited'' the same model and provided an approximation of its slow manifold while using Girimaji's technique of local reduction. Then, J. Vanneste \cite{24a} studied a ``conservative form'' of the Lorenz-Krishnamurthy model and developed distinct methods to derive the leading-order asymptotics of the late coefficients in the power-series expansion used for the construction of its slow manifold. But, to our knowledge, only the first-order and third-order approximations of this \textit{slow manifold} have been analytically obtained. So, the aim of this paper is to show that while using the \textit{Flow Curvature Method} one can provide the eighteenth-order approximation of the \textit{slow manifold} of the generalized Lorenz-Krishnamurthy model and the thirteenth-order approximation of the ``conservative'' Lorenz-Krishnamurthy model. The invariance of each \textit{slow manifold} is then established according to \textit{Darboux invariance theorem}.\\

This paper is organize as follows. The classical definitions of \textit{singularly perturbed systems} are briefly recalled in Sec. 2. Foundations of the \textit{Flow Curvature Method} are summarized in Sec. 3. More particularly, the definition of the \textit{flow curvature manifold} which provides an approximation of the \textit{slow manifold} associated with \textit{singularly perturbed systems} is presented in Prop. 1. Then, invariance of the \textit{flow curvature manifold} is stated according to Darboux theorem and to Prop. 2. In Sec. 4, application of the \textit{Flow Curvature Method} enables to provide the eighteenth-order approximation of the \textit{slow manifold} associated with the \textit{generalized Lorenz-Krishnamurthy model} \cite{2a} and the thirteenth-order approximation of the \textit{conservative Lorenz-Krishnamurthy model} \cite{24a}. The invariance of each \textit{slow manifold} is then established according to \textit{Darboux invariance theorem} and to Prop. 2.

\newpage

\section{Singularly perturbed dynamical systems}

Following the approach of C.K.R.T. Jones \cite{13} and Kaper
\cite{14} some fundamental concepts and definitions for systems of
ordinary differential equations with two time scales, i.e., for
\textit{singularly perturbed dynamical systems} are briefly recalled.

In the following we consider a dynamical systems theory for systems
of differential equations of the form:

\begin{equation}
\label{eq1}
\left\{ {{\begin{array}{*{20}c}
 {{\vec {x}}' = \vec {f}\left( {\vec {x},\vec {z},\varepsilon }
\right)\mbox{ }} \hfill \\
 {{\vec {z}}' = \varepsilon \vec {g}\left( {\vec {x},\vec {z},\varepsilon }
\right)} \hfill \\
\end{array} }} \right.
\end{equation}

where $\vec {x} \in \mathbb{R}^m$, $\vec {z} \in \mathbb{R}^p$,
$\varepsilon \in \mathbb{R}^ + $, and the prime denotes
differentiation with respect to the independent variable $t$. The
functions $\vec {f}$ and $\vec {g}$ are assumed to be $C^\infty$
functions\footnote{In certain applications these functions will be
supposed to be $C^r$, $r \geqslant 1$} of $\vec {x}$,
$\vec {z}$ and $\varepsilon$ in $U\times I$, where $U$ is an open
subset of $\mathbb{R}^m\times \mathbb{R}^p$ and $I$ is an open
interval containing $\varepsilon = 0$.

In the case when $\varepsilon \ll 1$, i.e., is a small positive
number, the variable $\vec {x}$ is called \textit{fast} variable,
and $\vec {z}$ is called \textit{slow} variable. Using Landau's
notation: $O ( \varepsilon ^l )$ represents a real
polynomial in $\varepsilon $ of $l$ degree, with $l \in \mathbb{Z}$,
it is used to consider that generally $\vec {x}$ evolves at an
$O\left( 1 \right)$ rate; while $\vec {z}$ evolves at an $O\left(
\varepsilon \right)$ \textit{slow} rate.

Reformulating the system (\ref{eq1}) in terms of the rescaled
variable $\tau = \varepsilon t$, we obtain:

\begin{equation}
\label{eq2}
\left\{ {{\begin{array}{*{20}c}
 {\varepsilon \dot {\vec {x}} = \vec {f}\left( {\vec {x},\vec{z},\varepsilon } \right)\mbox{ }} \hfill \\
 {\dot {\vec {z}} = \vec {g}\left( {\vec {x},\vec {z},\varepsilon } \right)}
\hfill \\
\end{array} }} \right.
\end{equation}

The dot $\left( \cdot \right)$ represents the derivative with
respect to the new independent variable $\tau $.

The independent variables $t$ and $\tau $ are referred to the
\textit{fast} and \textit{slow} times, respectively, and (\ref{eq1})
and (\ref{eq2}) are called the \textit{fast} and \textit{slow}
systems, respectively. These systems are equivalent whenever
$\varepsilon \ne 0$, and they are labeled \textit{singular
perturbation problems} when $\varepsilon \ll 1$, i.e., is a small
positive parameter. The label ``singular'' stems in part from the
discontinuous limiting behavior in the system (\ref{eq1}) as
$\varepsilon \to 0$.

In such case, the system (\ref{eq1}) reduces to an $m$-dimensional
system called \textit{reduced fast system}, with the variable $\vec
{z}$ as a constant parameter:

\begin{equation}
\label{eq3}
\left\{ {{\begin{array}{*{20}c}
 {{\vec {x}}' = \vec {f}\left( {\vec {x}, \vec {z}, 0 } \right)}
\hfill \\
 {{\vec {z}}' = \vec {0}\mbox{ }} \hfill \\
\end{array} }} \right.
\end{equation}

System (\ref{eq2}) leads to the following differential-algebraic
system called \textit{reduced slow system} which dimension decreases
from $m + p$ to $p$:

\begin{equation}
\label{eq4}
\left\{ {{\begin{array}{*{20}c}
 {\vec {0} = \vec {f}\left( {\vec {x},\vec {z},0 } \right)} \hfill
\\
 {\dot {\vec {z}} = \vec {g}\left( {\vec {x},\vec {z},0 } \right)}
\hfill \\
\end{array} }} \right.
\end{equation}

\newpage

By exploiting the decomposition into \textit{fast} and \textit{slow}
reduced systems (\ref{eq3}) and (\ref{eq4}), the geometric approach
reduced the full \textit{singularly perturbed system} to separate
lower-dimensional regular perturbation problems in the \textit{fast}
and \textit{slow} regimes, respectively.

\section{Fenichel geometric theory}

Fenichel geometric theory for general systems (\ref{eq1}), i.e., a theorem providing
conditions under which \textit{normally hyperbolic invariant manifolds} in system (\ref{eq1}) persist when the perturbation is
turned on, i.e., when $0 < \varepsilon \ll 1$ is briefly recalled in this subsection. This theorem concerns only compact manifolds with boundary.

\subsection{Normally hyperbolic manifolds}

Let's make the following assumptions about system (\ref{eq1}):

\textbf{(H}$_{1}$\textbf{)} \textit{The functions }$\vec
{f}$ \textit{and} $\vec {g}$ \textit{are} $C^\infty $\textit{ in
}$U\times I$\textit{, where U is an open subset of
}$\mathbb{R}^m\times \mathbb{R}^p$\textit{ and I is an open interval
containing }$\varepsilon = 0.$

\textbf{(H}$_{2}$\textbf{)} \textit{There exists a set }$M_0
$\textit{ that is contained in }$\left\{ {\left( {\vec {x},\vec {y}}
\right):\vec {f}\left( {\vec {x},\vec {y},0} \right) = {\vec {0}}}
\right\}$\textit{ such that }$M_0 $\textit{ is a compact manifold
with boundary and }$M_0 $\textit{ is given by the graph of a
}$C^\infty $\textit{ function }$\vec {y} = \vec {Y}_0 \left( \vec
{x} \right)$\textit{ for }$\vec {x} \in D$\textit{, where }$D
\subseteq \mathbb{R}^p$\textit{ is a compact, simply connected
domain and the boundary of D is an} $\left( {p - 1} \right)$
\textit{dimensional} $C^\infty$ \textit{ submanifold. Finally, the
set D is overflowing invariant with respect to (\ref{eq2}) when
}$\varepsilon = 0. $

\newpage

\textbf{(H}$_{3}$\textbf{)} $M_0 $\textit{ is normally hyperbolic
relative to (\ref{eq3}) and in particular it is required for all
points }$\vec {p} \in M_0 $\textit{, that there are k (resp. l)
eigenvalues of }$D_{\vec{y}} \vec{f}\left( {\vec{p},0}
\right)$\textit{ with positive (resp. negative) real parts bounded
away from zero, where }$k + l = m.$

\subsection{Fenichel persistence theory for singularly perturbed
systems}

For compact manifolds with boundary, Fenichel's persistence theory
states that, provided the hypotheses $(H_{1})-(H_{3})$ are
satisfied, the system (\ref{eq1}) has a \textit{slow} (or center)
\textit{manifold}, and this \textit{slow manifold} has \textit{fast}
stable and unstable \textit{manifolds}.

\textbf{\textit{Theorem for compact manifolds with boundary:}}

\textit{Let system (\ref{eq1}) satisfy the conditions
(H}$_{1})-(H_{3}$\textit{). If }$\varepsilon > 0$\textit{ is
sufficiently small, then there exists a function }$\vec {Y}\left(
{\vec {x},\varepsilon } \right)$\textit{ defined on D such that the
manifold }$M_\varepsilon = \left\{ {\left( {\vec {x},\vec {y}}
\right):\vec {y} = \vec {Y}\left( {\vec {x},\varepsilon } \right)}
\right\}$ \textit{is locally invariant under (\ref{eq1}). Moreover,
}$\vec {Y}\left( {\vec {x},\varepsilon } \right)$\textit{ is
}$C^r$\textit{ for any }$r < + \infty $, \textit{ and
}$M_\varepsilon $\textit{ is }$C^rO\left( \varepsilon \right)$
\textit{close to }$M_0 $\textit{. In addition, there exist perturbed
local stable and unstable manifolds of }$M_\varepsilon $\textit{.
They are unions of invariant families of stable and unstable fibers
of dimensions l and k, respectively, and they are} $C^rO\left(
\varepsilon \right)$ \textit{close for all} $r < + \infty$ \textit{,
to their counterparts.}

\textit{\textbf{Proof.}}\hspace{0.1in} For proof of this theorem see
Fenichel \cite{{5},{6},{7},{8}}.

The label slow manifold is attached to $M_\varepsilon $ because the
magnitude of the vector field restricted to $M_\varepsilon $ is
$O\left( \varepsilon \right)$, in terms of the fast independent
variable t.

So persistent manifolds are labeled \textit{slow manifolds}, and the
proof of their persistence is carried out by demonstrating that the
local stable and unstable manifolds of $M_0$ also persist as locally
invariant manifolds in the perturbed system, i.e., that the local
hyperbolic structure persists, and then the slow manifold is
immediately at hand as a locally invariant manifold in the
transverse intersection of these persistent local stable and
unstable manifolds.

\subsection{Geometric singular perturbation theory}

Earliest geometric approaches to \textit{singularly perturbed
systems} have been developed by Cole \cite{2}, O'Malley
\cite{{16},{17}}, Fenichel \cite{{5},{6},{7},{8}} for the
determination of the \textit{slow manifold} equation.

Generally, Fenichel theory enables to turn the problem for
explicitly finding functions $\vec {y} = \vec {Y}\left( {\vec
{x},\varepsilon } \right)$ whose graphs are locally invariant
\textit{slow manifolds} $M_\varepsilon $ of system (\ref{eq1}) into
regular perturbation problem \cite[p. 112]{16}. Invariance of the manifold
$M_\varepsilon $ implies that $\vec {Y}\left( {\vec {x},\varepsilon
} \right)$ satisfies:

\begin{equation}
\label{eq5}
D_{\vec {x}} \vec {Y}\left( {\vec
{x},\varepsilon } \right)\vec {f}\left(\vec{x}, {\vec {Y}\left( {\vec
{x},\varepsilon } \right),\varepsilon } \right) = \varepsilon  \vec
{g}\left(\vec{x}, {\vec {Y}\left( {\vec {x},\varepsilon } \right),\varepsilon } \right)
\end{equation}

According to Guckenheimer \textit{et al.} \cite[p. 131]{10b}, this (partial)
differential equation for $\vec {Y}\left( {\vec {x}, \varepsilon} \right)$ cannot be
solved exactly. So, its solution can be approximated arbitrarily closely as
a Taylor series at $\left( {\vec {x}, \varepsilon} \right) = \left( \vec {0}, 0 \right)$.

Then, the following perturbation expansion is plugged:

\begin{equation}
\label{eq6}
\vec{Y}\left( {\vec {x},\varepsilon } \right) = \vec {Y}_0 \left( \vec
{x} \right) + \varepsilon \vec {Y}_1 \left( \vec {x} \right) + O\left( {\varepsilon^2} \right)
\end{equation}

into (\ref{eq5}) to solve order by order for $\vec {Y}\left( {\vec {x},\varepsilon } \right)$. The
Taylor series expansion \cite{17} for $\vec {f}\left(\vec{x}, {\vec {Y}\left(
{\vec {x},\varepsilon } \right),\varepsilon } \right)$ and $\vec {g}\left(\vec{x}, {\vec {Y}\left(
{\vec {x},\varepsilon } \right),\varepsilon } \right)$ up
to terms of order two in $\varepsilon $ reads:

\[
\begin{aligned}
& \vec {f}\left(\vec{x}, {\vec {Y}\left( {\vec {x},\varepsilon } \right),\varepsilon } \right)  =  \vec {f}\left(\vec{x}, {\vec {Y}_0 \left( \vec {x} \right),0} \right) +  \varepsilon \left( D_{\vec {y}} \vec {f}\left(\vec{x}, {\vec {Y}_0 \left( \vec {x} \right),0} \right) \vec {Y}_1\left( \vec {x} \right) + \frac{\partial \vec {f}}{\partial \varepsilon }\left(\vec{x}, {\vec {Y}_0 \left( \vec {x} \right),0} \right) \right) \\
& \vec {g}\left(\vec{x}, {\vec {Y}\left( {\vec {x},\varepsilon } \right),\varepsilon } \right)  =  \vec {g}\left(\vec{x}, {\vec {Y}_0 \left( \vec {x} \right),0} \right) +  \varepsilon \left( D_{\vec {y}} \vec {g}\left(\vec{x}, {\vec {Y}_0 \left( \vec {x} \right),0} \right) \vec {Y}_1\left( \vec {x} \right) + \frac{\partial \vec {g}}{\partial \varepsilon }\left(\vec{x}, {\vec {Y}_0 \left( \vec {x} \right),0} \right) \right)
\end{aligned}
\]

\textbullet \hspace{0.1in} At order $\varepsilon ^0$, Eq. (\ref{eq5}) gives:

\begin{equation}
\label{eq7}
D_{\vec {x}} \vec{Y_0}\left( \vec{x}\right)  \vec {f}(\vec{x}, {\vec {Y}_0 \left( \vec {x} \right),0}) = \vec {0}
\end{equation}

which defines $\vec {Y}_0 \left( \vec {x} \right)$ due to the invertibility of $D_{\vec {y}} \vec {f}$ and the \textit{Implicit Function Theorem}.

\textbullet \hspace{0.1in} The next order $\varepsilon ^1$ provides:

\begin{equation}
\label{eq8}
D_{\vec {x}} \vec {Y}_0 \left( \vec {x} \right) \left[ D_{\vec {y}} \vec {f} (\vec{x}, {\vec {Y}_0 \left( \vec {x} \right),0} ) \vec{Y_1}\left( \vec{x} \right) + \frac{\partial \vec {f}}{\partial \varepsilon } \right] = \vec {g}(\vec{x}, {\vec {Y}_0 \left( \vec {x} \right),0} )
\end{equation}

which yields $\vec {Y}_1 \left( \vec {x} \right)$ and so forth.

So, regular perturbation theory makes it possible to build an approximation of locally
invariant \textit{slow manifolds} $M_\varepsilon$. Thus, in the
framework of the \textit{Geometric Singular Perturbation Method},
three conditions are needed to characterize the \textit{slow
manifold} associated with \textit{singularly perturbed system}:
existence, local invariance and determination. Existence and local
invariance of the \textit{slow manifold} are stated according to
Fenichel theorem for compact manifolds with boundary while
\textit{asymptotic expansions} provide its equation up to the order
of the expansion.

\newpage

\section{Flow Curvature Method}

In this section, one of the main results of the \textit{Flow Curvature Method} and based on the use of local properties of \textit{curvatures} inherent to \textit{Differential Geometry} is briefly presented (for more details see \cite{9,9d}). According to this method, the highest \textit{curvature of the flow}, i.e. the $(n - 1)^{th}$ \textit{curvature} of \textit{trajectory curve} integral of $n$--dimensional \textit{singularly perturbed dynamical systems} (\ref{eq1}) defines a $(n - 1)$--dimensional \textit{manifold} associated with this system and called \textit{flow curvature manifold}. We have the following result:

\subsection{Slow manifold equation}

\begin{prop}
\label{prop1}
The location of the points where the $(n - 1)^{th}$ \textit{curvature of the flow}, i.e. the \textit{curvature of the
trajectory curve} $\vec {X}$, integral of any $n$--dimensional singularly perturbed dynamical systems (\ref{eq1}) vanishes, provides a $k$--order approximation in $\varepsilon$ of its \textit{slow manifold} $M_{\varepsilon}$ the equation of which reads

\begin{equation}
\label{eq9} \phi ( {\vec {X}, \varepsilon} ) = \dot { \vec
{X} } \cdot ( {\ddot {\vec {X}} \wedge \dddot{\vec {X}}\wedge \ldots
\wedge \mathop {\vec {X}}\limits^{\left( n \right)} } ) = \det( {
\dot {\vec {X}},\ddot {\vec {X}}, \dddot{\vec {X}},\ldots ,\mathop
{\vec {X}}\limits^{\left( n \right)} } ) = 0
\end{equation}

where $\mathop {\vec {X}}\limits^{\left( n \right)}$ represents the
time derivatives up to order $n$ of $\vec {X}  = (\vec{x}, \vec{z})^t$.
\end{prop}

\textit{\textbf{Note.}}\hspace{0.1in}

$k$-order approximation depends on the number of $\varepsilon$ contained in the vector field. We will see in Sec. 4 that for the Lorenz-Krishnamurthy model $k=18$.

\newpage

\textit{\textbf{Proof.}}\hspace{0.1in}

While the \textit{slow invariant manifold} analytical equation (\ref{eq6}) given by the \textit{Geometric Singular Perturbation Method} is an \textit{explicit} equation, the \textit{slow invariant manifold} analytical equation (\ref{eq9}) obtained according to the \textit{Flow Curvature Method} is an \textit{implicit} equation. So, in order to compare the latter with the former it is necessary to plug the following perturbation expansion: $\vec {Y}\left( {\vec {x},\varepsilon } \right) = \vec {Y}_0 \left( \vec {x} \right) + \varepsilon \vec {Y}_1 \left( \vec {x} \right) + O\left( {\varepsilon ^2} \right)$ into (\ref{eq9}). Thus, solving order by order for $\vec {Y}\left( {\vec {x},\varepsilon } \right)$ will transform (\ref{eq9}) into an \textit{explicit analytical equation} enabling the comparison with (\ref{eq6}). The Taylor series expansion for $\phi (\vec{X}, \varepsilon ) = \phi (\vec{x}, \vec {Y}\left( {\vec {x},\varepsilon } \right), \varepsilon )$ up to terms of order one in $\varepsilon $ reads:

\begin{equation}
\label{eq10} \phi (\vec{X}, \varepsilon ) = \phi (\vec{x}, \vec {Y}_0
\left( \vec {x} \right),0 ) + \varepsilon D_{\vec{y}} \phi (\vec{x}, \vec {Y}_0 \left( \vec {x} \right),0 )\vec {Y}_1 \left( \vec {x} \right) + \varepsilon \frac{\partial \phi }{\partial \varepsilon } (\vec{x}, \vec {Y}_0 \left( \vec {x} \right),0 )
\end{equation}

\textbullet \hspace{0.1in} At order $\varepsilon ^0$, Eq. (\ref{eq10}) gives:

\begin{equation}
\label{eq11} \phi \left(\vec {x}, \vec {Y}_0 \left( \vec {x} \right),0 \right) = 0
\end{equation}

which defines $\vec {Y}_0 \left( \vec {x} \right)$ due to the
invertibility of $D_{\vec {y}} \phi $ and application of the
\textit{Implicit Function Theorem}.

\textbullet \hspace{0.1in} The next order $\varepsilon ^1$, provides:

\begin{equation}
\label{eq12} D_{\vec{y}} \phi (\vec{x}, \vec {Y}_0 \left( \vec {x} \right),0 )\vec {Y}_1 \left( \vec {x} \right) + \frac{\partial \phi }{\partial \varepsilon } (\vec{x}, \vec {Y}_0 \left( \vec {x} \right),0 ) = \vec {0}
\end{equation}

which yields $\vec {Y}_1 \left( \vec {x} \right)$ and so forth.

In order to prove that this equation is completely identical to Eq. (\ref{eq8}), let's rewrite it as follows:

\[
\vec {Y}_1 \left( \vec {x} \right) = - \left[ D_{\vec{y}} \phi (\vec{x}, \vec {Y}_0 \left( \vec {x} \right),0 ) \right]^{-1} \frac{\partial \phi }{\partial \varepsilon } (\vec{x}, \vec {Y}_0 \left( \vec {x} \right),0 )
\]

By application of the \textit{chain rule}, i.e., the derivative of $\phi (\vec{x}, \vec {Y}_0 \left( \vec {x} \right),0)$ with respect to the variable $\vec {y}$ and then with respect to $\varepsilon$, it can be stated that:

\[
\vec {Y}_1 \left( \vec {x} \right)  = - \left[ (D_{\vec {x}} \vec {f}) (D_{\vec {y}}\vec {f}) \right]^{-1} (D_{\vec {y}}\vec {f}) \vec {g}(\vec{x}, \vec {Y}_0 ( \vec {x} ),0 ) - \left[ D_{\vec {y}}\vec {f} \right]^{-1} D_{\varepsilon}\vec {f}(\vec{x}, \vec {Y}_0 \left( \vec {x} \right),0 )
\]

But, according to the \textit{Implicit Function Theorem} we have:

\[
(D_{\vec {x}} \vec {f}) = - (D_{\vec {y}} \vec {f}) (D_{\vec {x}} \vec {y}) = - (D_{\vec {y}} \vec {f}) (D_{\vec {x}} \vec {Y}_0 ( \vec {x} ))
\]

Then, by replacing into the previous equation we find:

\[
\vec {Y}_1 \left( \vec {x} \right)  = \left[ (D_{\vec {y}} \vec {f}) (D_{\vec {x}} \vec {Y}_0 ( \vec {x} )) (D_{\vec {y}}\vec {f}) \right]^{-1} (D_{\vec {y}}\vec {f}) \vec {g}(\vec{x}, \vec {Y}_0 ( \vec {x} ),0 ) - \left[ D_{\vec {y}}\vec {f} \right]^{-1} D_{\varepsilon}\vec {f}(\vec{x}, \vec {Y}_0 \left( \vec {x} \right),0 )
\]

After simplifications, we have:

\[
\vec{Y_1}\left( \vec{x} \right) = \left[ D_{\vec {x}} \vec {Y}_0 \left( \vec {x} \right) D_{\vec {y}} \vec {f} (\vec{x}, {\vec {Y}_0 \left( \vec {x} \right),0} ) \right]^{-1} \vec {g}(\vec{x}, {\vec {Y}_0 \left( \vec {x} \right),0} ) - \left[ D_{\vec {y}} \vec {f} (\vec{x}, {\vec {Y}_0 \left( \vec {x} \right),0} ) \right]^{-1} \frac{\partial \vec {f}}{\partial \varepsilon }
\]

\newpage

Finally, Eq. (\ref{eq12}) may be written as:

\[
D_{\vec {x}} \vec {Y}_0 \left( \vec {x} \right) \left[ D_{\vec {y}} \vec {f} (\vec{x}, {\vec {Y}_0 \left( \vec {x} \right),0} ) \vec{Y_1}\left( \vec{x} \right) + \frac{\partial \vec {f}}{\partial \varepsilon } \right] = \vec {g}(\vec{x}, {\vec {Y}_0 \left( \vec {x} \right),0} )
\]

Thus, identity between the ``\textit{slow manifold}'' equation given by the \textit{Geometric Singular Perturbation Method} and by the \textit{Flow Curvature Method} is proved up to first order term in $\varepsilon$.

\textit{\textbf{Note.}}\hspace{0.1in}
Let's notice that the \textit{slow invariant manifold} equation (\ref{eq9}) associated with $n$--dimensional \textit{singularly
perturbed systems} defined by the \textit{Flow Curvature Method} is a tensor of order $n$. As a consequence, it can only provide an approximation of  $n$--order in $\varepsilon$ of the \textit{slow invariant manifold} equation (\ref{eq6}).
Nevertheless, it is easy to show that the Lie derivative of the ``\textit{slow manifold}'' equation (\ref{eq9}) obtained by the \textit{Flow Curvature Method} can be written as:

\begin{equation}
\label{eq13}
L_{\vec {X}} \phi ( {\vec {X}, \varepsilon} ) = \dot { \vec {X} } \cdot ( {\ddot {\vec {X}} \wedge \dddot{\vec {X}}\wedge \ldots
\wedge \mathop {\vec {X}}\limits^{\left( n+1 \right)} } ) = \det( {\dot {\vec {X}},\ddot {\vec {X}}, \dddot{\vec {X}},\ldots ,\mathop
{\vec {X}}\limits^{\left( n+1 \right)} } ) = 0
\end{equation}

where $\mathop {\vec {X}}\limits^{\left( n + 1 \right)}$ represents the time derivatives up to order $(n + 1)$ of $\vec {X}  = (\vec{x}, \vec{y})^t$. So, Eq. (\ref{eq13}) defines a tensor of order $n+1$ which provides an approximation of $(n + 1)$-order in $\varepsilon$ of the \textit{slow invariant manifold} equation (\ref{eq6}). Thus, by taking the successive Lie derivatives of the ``\textit{slow manifold}'' equation (\ref{eq9}) we improve the order of the approximation up to an order corresponding to that of the Lie derivative.
As an example, according to Prop. 1, the ``\textit{slow manifold}'' equation of a two-dimensional singularly perturbed dynamical system reads:

\[
\phi(\vec{X}, \varepsilon) = \det(\dot {\vec {X}},\ddot {\vec {X}}) = 0
\]

where $\vec{X} = (x,y)^t$. This second-order tensor only provides a first order approximation in $\varepsilon$ of the \textit{slow invariant manifold} equation (\ref{eq6}). While its Lie derivative

\[
L_{\vec {X}} \phi(\vec{X}, \varepsilon) = \det(\dot {\vec {X}},\dddot {\vec {X}}) = 0
\]

which is third-order tensor gives a second-order approximation in $\varepsilon$. Thus, by applying a mathematical induction, the proof above can be extended to high-order approximations in $\varepsilon$.

\subsection{Invariance of the slow manifold}

The local invariance of the \textit{slow manifold} analytical equation defined by \textit{Flow Curvature Method} may be
stated while using the \textit{Tangent Linear System Approximation (T.L.S.A.)} associated with \textit{Darboux Invariance Theorem}. \textit{Tangent Linear System Approximation} has been introduced by Rossetto \textit{et al.} \cite{21} in order to compute
the \textit{slow manifold} equation of \textit{singularly perturbed systems}. This approximation consists in replacing, in the vicinity of the \textit{singular approximation}, the \textit{singularly perturbed system} by the corresponding \textit{tangent linear system}. \textit{Tangent
Linear System Approximation} may thus be viewed as a local invariance condition of the \textit{slow manifold}.

\subsubsection{Tangent linear system approximation (T.L.S.A.)}

The tangent linear system approximation (T.L.S.A.) states that, in
the vicinity of the \textit{singular approximation} associated with \textit{singularly perturbed
systems} (\ref{eq1}), the functional jacobian matrix of such systems is locally stationary, i.e.,

\begin{equation}
\label{eq14} \frac{dJ}{dt} = 0
\end{equation}

\subsubsection{Lie Derivative - Darboux Invariance Theorem}

Let $\phi $ a $C^1$ function defined in a compact E included in
$\mathbb{R}$ and $\vec {X}\left( t \right)$ the integral of the
dynamical system defined by (\ref{eq1}). The Lie derivative is
defined as follows:

\begin{equation}
\label{eq15} L_{\vec {X}} \phi = \dot{\vec{ X }} \cdot
\overrightarrow \nabla\phi = \sum\limits_{i = 1}^n {\frac{\partial
\phi }{\partial x_i }\dot {x}_i } = \frac{d\phi }{dt}
\end{equation}

\textbf{\textit{Darboux Invariance Theorem:}}

An \textit{invariant manifold} is defined by $\phi ( \vec {X}, \varepsilon ) = 0$ where $\phi $ is a $C^1$ in an
open set U and such there exists a $C^1$ function denoted $k( \vec {X} )$ and called cofactor which satisfies

\begin{equation}
\label{eq16} L_{\vec {X}} \phi ( \vec {X}, \varepsilon ) =
k ( \vec {X} ) \phi ( \vec {X}, \varepsilon ) \quad \mbox{for all} \quad \vec {X} \in U
\end{equation}

\textit{\textbf{Proof.}}\hspace{0.1in}

The proof of this theorem is in Darboux \cite{3}. However, let's prove that both Darboux and Fenichel's invariance are exactly identical.

According to \textit{Fenichel's persistence theorem} a \textit{slow invariant manifold} $M_{\varepsilon}$ may be written as an \textit{explicit} function: $\vec{y} = \vec {Y}\left( {\vec {x},\varepsilon } \right)$ the invariance of which implies that $\vec {Y}\left( {\vec {x},\varepsilon } \right)$ satisfies:

\begin{equation}
\label{eq17}
D_{\vec {x}} \vec {Y}\left( {\vec {x},\varepsilon } \right)\vec{f}(\vec{x}, {\vec {Y}\left( {\vec {x},\varepsilon } \right),\varepsilon } ) = \varepsilon \vec {g}(\vec{x}, {\vec {Y}\left( {\vec {x},\varepsilon } \right),\varepsilon } )
\end{equation}

Let's write the \textit{slow manifold} $M_{\varepsilon}$ as an \textit{implicit function} by posing:

\begin{equation}
\label{eq18}
\phi(\vec {x},\vec {y},\varepsilon ) = \vec{y} - \vec {Y}\left( \vec {x},\varepsilon  \right)
\end{equation}

According to \textit{Darboux Invariance Theorem} $M_{\varepsilon}$ is invariant if and only if:

\begin{equation}
\label{eq19}
L_{\vec V } \phi ( \vec {x},\vec {y},\varepsilon  ) = k(\vec {x},\vec {y},\varepsilon ) \phi ( \vec {x},\vec {y},\varepsilon )
\end{equation}

Plugging Eq. (\ref{eq18}) into the Lie derivative (\ref{eq19}) leads to:

\[
L_{\vec V } \phi ( \vec {x},\vec {y},\varepsilon ) = \dot{\vec{y}} - D_{\vec {x}}\vec {Y}\left( \vec {x},\varepsilon \right) \dot{\vec{x}} = k(\vec {x},\vec {y},\varepsilon ) \phi ( \vec {x},\vec {y},\varepsilon )
\]

which may be written according to Eq. (\ref{eq1}):

\[
L_{\vec V } \phi ( \vec {x},\vec {y},\varepsilon ) = \varepsilon \vec{g}(\vec{x},\vec{y},\varepsilon) - D_{\vec {x}}\vec {Y}\left( {\vec {x},\varepsilon } \right) \vec{f}(\vec{x},\vec{y},\varepsilon) = k(\vec {x},\vec {y},\varepsilon) \phi ( \vec {x},\vec {y},\varepsilon)
\]

\vspace{0.1in}
Evaluating this Lie derivative in the location of the points where $\phi(\vec {x},\vec {y},\varepsilon ) = 0$, i.e. $\vec{y} = \vec {Y}\left( \vec {x},\varepsilon  \right)$ leads to:

\[
L_{\vec V } \phi (\vec{x}, \vec {Y}\left( \vec {x},\varepsilon \right),\varepsilon ) =  \varepsilon \vec{g}(\vec{x}, \vec {Y}\left( \vec {x},\varepsilon  \right),\varepsilon) - D_{\vec {x}}\vec {Y}\left( {\vec {x},\varepsilon } \right) \vec{f}(\vec{x},\vec {Y}\left( \vec {x},\varepsilon  \right),\varepsilon) = 0
\]

which is exactly identical to Eq. (\ref{eq5}) proposed by Fenichel\index{Fenichel}.

Now, let's prove the invariance of the \textit{flow curvature manifold} (\ref{eq9}), \textit{i.e.}, the invariance of the ``\textit{slow manifold equation}'' defined by Prop. 1.

\begin{prop}
\label{prop2}
\textit{The flow curvature manifold defined by }$\phi ( {\vec {X}} )=0$\textit{ where }$\phi $\textit{ is a }$C^1$\textit{ in an open set U is invariant with respect to the flow of (\ref{eq1}) if there exists a }$C^1$\textit{ function denoted }$k ( {\vec {X}} )$ \textit{and called cofactor which satisfies:}

\begin{equation}
\label{eq20}
L_{\vec V } \phi ( {\vec {X}} )= k ( {\vec {X}} )\phi ( {\vec {X}} )
\end{equation}

\vspace{0.1in}

\textit{for all} $\vec {X}\in U$ \textit{and where} $L_{\vec V} \phi =\vec V \cdot \vec \nabla \phi =\sum\limits_{i=1}^n {\dfrac{\partial \phi}{\partial x_i }\dot {x}_i } =\dfrac{d\phi }{dt}.$

\end{prop}

\textit{\textbf{Proof.}}\hspace{0.1in} Lie derivative of the \textit{flow curvature manifold} (\ref{eq9}) reads:

\begin{equation}
\label{eq21}
L_{\vec V } \phi ( {\vec {X}} )=\dot {\vec {X}}\cdot ( {\ddot {\vec {X}} \wedge \dddot{\vec {X}}\wedge \ldots \wedge \mathop {\vec {X}}\limits^{\left( {n+1} \right)} } )
\end{equation}

From the identity $\ddot {\vec {X}}=J\dot {\vec {X}}$ where $J$ is the functional jacobian matrix associated with any $n$--dimensional \textit{singularly perturbed system} (\ref{eq1}) we find that:

\begin{equation}
\label{eq22}
\mathop {\vec {X}}\limits^{\left( {n+1} \right)} =J^n\dot {\vec {X}} \mbox{\quad if \quad} \frac{dJ}{dt}=0
\end{equation}

where $J^n$ represents the $n^{th}$ power of $J$, e.g., $\ddot {\vec
{X}} = J\dot {\vec {X}}$, $\dddot {\vec {X}} = J\ddot {\vec {X}}$, \ldots

Then, it follows that:

\begin{equation}
\label{eq23}
\mathop {\vec {X}}\limits^{\left( {n+1} \right)} =J\mathop {\vec {X}}\limits^{\left( n \right)}
\end{equation}

Replacing $\mathop {\vec {X}}\limits^{\left( {n+1} \right)} $ in Eq. (\ref{eq14}) with Eq. (\ref{eq16}) we have:

\begin{equation}
\label{eq24}
L_{\vec V } \phi ( {\vec {X}} ) = \dot {\vec {X}}\cdot ( {\ddot {\vec {X}}\wedge \dddot{\vec {X}}\wedge \ldots \wedge J\mathop {\vec{X}}\limits^{\left( n \right)} } )
\end{equation}

The right hand side of this Eq. (\ref{eq24}) can be written:

\[
J \dot {\vec {X}}\cdot ( {\ddot {\vec {X}}\wedge \dddot{\vec {X}}\wedge \ldots \wedge \mathop {\vec{X}}\limits^{\left( n \right)} } ) + \dot {\vec {X}}\cdot ( J {\ddot {\vec {X}}\wedge \dddot{\vec {X}}\wedge \ldots \wedge \mathop {\vec{X}}\limits^{\left( n \right)} } ) + \ldots + \dot {\vec {X}}\cdot ( {\ddot {\vec {X}}\wedge \dddot{\vec {X}}\wedge \ldots \wedge J\mathop {\vec{X}}\limits^{\left( n \right)} } )
\]

According to Eq. (\ref{eq23}) all terms are null except the last one. So, by taking into account identity (\ref{eqA1}) established in Appendix we find:

\[
L_{\vec V } \phi ( {\vec {X}} ) = Tr\left( J \right) \dot {\vec {X}}\cdot ( {\ddot {\vec {X}}\wedge \dddot{\vec {X}}\wedge \ldots
\wedge \mathop {\vec {X}}\limits^{\left( n \right)} } ) = Tr\left( J \right)\phi ( {\vec {X}} ) = k( {\vec {X}} ) \phi ( {\vec {X}} )
\]

where $k( {\vec {X}} )=Tr\left( J \right)$ represents the trace of the functional jacobian matrix.

So, according to Prop. \ref{prop2} invariance of the \textit{slow manifold} analytical equation of any $n$--dimensional \textit{singularly perturbed dynamical system} is established provided that the functional jacobian matrix is locally stationary.

\section{The Lorenz-Krishnamurthy Slow Invariant Manifold}

Let's consider the model introduced by the late E.N. Lorenz \cite{15a} and usually referred
to as the Lorenz five-mode model \cite{1a} or as the Lorenz-Krishnamurthy model \cite{{15a},{15b},{15c}}. This model, obtained by truncation of the rotating shallow-water equations, governs the dynamics of a triad of vortical modes, with amplitudes $(u,v,w)$, coupled to a gravity mode described by $(x, y)$. In 1996, Camassa \textit{et al.} \cite{2a} proposed a \textit{generalized LK model} presented in the next section. Starting from this model we will provide the eighteenth-order approximation of its \textit{slow manifold} while using the \textit{Flow Curvature Method}, the invariance of which will be established according to Darboux theroem. Moreover, by posing $\delta = 0$ in our \textit{slow manifold} analytical equation we will find again the first-order approximation of the \textit{slow manifold} given by Camassa \textit{et al.} \cite[p. 3263]{2a}.

\subsection{The generalized LK model}

According to Camassa \textit{et al.} \cite{2a}, the Lorenz-Krishnamurthy model \cite{15b} can be written as:

\begin{equation}
\label{eq25}
\left\{ {\begin{array}{*{20}c}
\dot{x} =  -y - \kappa x\\
\dot{y} = x + \epsilon u v - \kappa y\\
\dot{u} = -v w + \epsilon v y - \alpha u \\
\dot{v} = u w - \epsilon u y - \alpha v + \alpha F\\
\dot{w} = - u v - \alpha w
\end{array} } \right.
\end{equation}

where $\epsilon \geqslant 0$ is the coupling parameter between the \textit{Rossby Wave} $(u,v,w)$ and \textit{Gravity Wave} $(x,y)$, $\alpha$ is a parameter introduced to model the dissipation and controls the damping of Rossby mode, while $\kappa$ controls the damping of the gravity mode. $F > 0$ represents the forcing parameter which is assumed to be much smaller than unity. By posing $\alpha = a$, $\epsilon = b$ and $\kappa = \alpha$ in Eqs. (\ref{eq25}) one finds again the Lorenz-Krishnamurthy model \cite{15b}. Then, by making the following variables changing:

\[
x \to \delta^2 x \mbox{, } y \to \delta^2 y \mbox{, } u \to \delta u \mbox{, } v \to \delta v \mbox{, } w \to \delta w
\]

and while posing for convenience $\alpha = \delta^2$, Camassa \textit{et al.} \cite{2a} obtained the following system:

\begin{equation}
\label{eq26}
\left\{ {\begin{array}{*{20}c}
\dot{x} =  -y - \kappa x\\
\dot{y} = x + \epsilon u v - \kappa y\\
\dot{u} = \delta ( -v w + \delta \epsilon v y - \delta u) \\
\dot{v} = \delta ( u w - \delta \epsilon u y - \delta v + F)\\
\dot{w} = -\delta ( u v + \delta w)
\end{array} } \right.
\end{equation}

\textit{\textbf{Note.}}\hspace{0.1in}

Reformulating the \textit{fast} system (\ref{eq26}) in terms of the rescaled variable $\tau = \delta t$ we obtain a \textit{two-time scales singularly perturbed dynamical system} of the same form as the \textit{slow} system (\ref{eq2}) for which variables $\vec{x} = (x,y)^t$ and $\vec{z} = (u,v,w)^t$ are respectively \textit{fast} and \textit{slow}.

\begin{equation}
\label{eq27}
\left\{ {\begin{array}{*{20}c}
\delta \dot{x} =  -y - \kappa x\\
\delta \dot{y} = x + \epsilon u v - \kappa y\\
\dot{u} =  -v w + \delta \epsilon v y - \delta u \\
\dot{v} =  u w - \delta \epsilon u y - \delta v + F \\
\dot{w} = - ( u v + \delta w)
\end{array} } \right.
\end{equation}

In the first article in which he introduced his five-mode model, Lorenz \cite[p. 1548]{15} wrote:

\begin{quote}
``In the trivial case where the Rossby and gravity waves are completely uncoupled, i.e., where the system of equations degenerates into two systems, one governing Rossby waves and one governing gravity waves, the slow manifold obviously exists and is obtained simply by equating all of the gravity-waves variables to zero.''
\end{quote}

Following this idea Camassa \textit{et al.} \cite[p. 3263]{2a} gave the zero-order approximation in $\delta$ of the \textit{slow manifold} associated with the system (\ref{eq27}):

\begin{equation}
\label{eq28}
\begin{aligned}
x & = - \epsilon \frac{ u v}{1 + \kappa^2} + O(\delta)  \\
y & = \kappa \epsilon \frac{u v}{1 + \kappa^2} + O(\delta) \\
\end{aligned}
\end{equation}

This \textit{slow manifold} parametrized as a graph over the $(u,v,w)$ space is also called \textit{singular approximation} (see Sec. 2 for definition) since it is obtained by posing $\delta = 0$ in the two first equations of system (\ref{eq27}).

\newpage

After the publication of the article entitled ``On the Nonexistence of a Slow Manifold'' in which Lorenz and Krishnamurthy \cite{15b} concluded that the slow manifold of such model ``does not exist'', Jacobs \cite{12a}, Boyd \cite{1a}, Fowler and Kember \cite{12a} and then Camassa and Tin \cite{2a} proved the existence of a \textit{slow manifold} in this model and gave approximations of its equation at first orders while using the ``singular perturbation scheme known as the method of multiple scales.\footnote{Boyd \cite[p. 1058]{1a}}''. Although these authors stated that a \textit{slow manifold} can be constructed via formal series, such a long and tedious asymptotic procedure of systematic identification order-by-order is expected to diverge as previously recalled. Recently, Sudheer \textit{et al.} \cite{17a} and Vanneste \cite{24a} have proposed alternative techniques for the construction of the \textit{slow manifold} of such model. We will show now that one can obtain high-orders approximation of this \textit{slow manifold} while using the \textit{Flow Curvature Method}.

Thus, according to Prop. 1 the \textit{slow manifold} equation (\ref{eq9}) associated with the \textit{generalized LK model} (\ref{eq27}) reads\footnote{See http://ginoux.univ-tln.fr for complete equation.}:

\begin{equation}
\label{eq29}
\phi(\vec{X}, \delta) = \det(\dot {\vec {X}},\ddot {\vec {X}},\dddot {\vec {X}},\ddddot {\vec {X}},\mathop{\vec {X}}\limits^{\ldots..} ) = 0
\end{equation}

where $\vec {X}  = (\vec{x}, \vec{z})^t$. Then, it can be verified that the time derivative of the functional jacobian matrix of system (\ref{eq26}) evaluated when  $\delta \to 0$ is a zero matrix. So, from \textit{Darboux Invariance Theorem} we can conclude that in the $\delta$-vicinity of the \textit{singular approximation} the \textit{slow manifold} is invariant.

\newpage

The implicit equation (\ref{eq29}) is a polynomial of degree $10$ for $u$, $v$ and $w$, of degree $5$ for $x$ and $11$ for $y$ and represents the eighteenth-order approximation in $\delta$ of the \textit{slow manifold} of the \textit{generalized LK model} (\ref{eq27}).

By posing $\delta = 0$ in the above Eq. (\ref{eq29}) we find that:

\begin{equation}
\label{eq30}
\phi(\vec{X}, 0) = v (u^2 - w^2) \left[ (x + \epsilon \frac{ u v}{1 + \kappa^2})^2 +(y - \kappa \epsilon \frac{u v}{1 + \kappa^2})^2 \right] = 0
\end{equation}

This equation is made of a product of \textit{invariant manifolds} as it is easy to verify according to \textit{Darboux Invariance Theorem}. Let's compute the Lie derivative of the first and second term, when $\delta \to 0$ we have:

\begin{equation}
\label{eq31}
\begin{aligned}
L_{\vec X } ( v ) & = 0 \\
L_{\vec X } ( u^2 - w^2 ) & = 0 \\
\end{aligned}
\end{equation}

Let's notice that the third term of Eq. (\ref{eq30}) is nothing else but the zero-order approximation in $\delta$ (\textit{singular approximation}) of the \textit{slow manifold} (see Eq. (\ref{eq28})) given by Camassa \textit{et al.} \cite[p. 3263]{2a} which is also invariant when $\delta \to 0$.

But, according to Leith \cite[p. 960]{15aaa}, the decomposition into \textit{fast} and \textit{slow} modes enables to define a three-dimensional submanifold of the state space parametrized by $(u,v,w)$ and that he called \textit{slow manifold}. So, let's pose $x \to 0$ and $y \to 0$ in the above Eq. (\ref{eq29}) we find that:

\begin{equation}
\label{eq32}
\begin{aligned}
\phi(u,v,w, \delta) & = u^2 w^2 (v^2 \delta^2 (1 + \epsilon^2) - (\delta^2 - \kappa)^2) - w^2 (u^2 + v^2 w^2 \delta^2) \\
& +  u^4 (1 + (\delta^2 - \kappa)^2) + u v w \delta (\delta^2 - \kappa) (2 w^2 - u^2 (2 + \epsilon^2)) = 0
\end{aligned}
\end{equation}

\newpage

In addition to the \textit{invariant manifolds} (\ref{eq31}) highlighted above we find another manifold. Let's compute its Lie derivative when $\delta \to 0$ we obtain:

\begin{equation}
\label{eq33}
L_{\vec X } \phi(u,v,w, 0)  = u^2 (u^2 - w^2)(1 + \kappa^2)
\end{equation}

Thus, we deduce that this manifold is \textit{locally} invariant, i.e. is invariant in the vicinity of the manifold defined by $u^2 - w^2  = 0$.

Now, by posing in system (\ref{eq27}) $\kappa = \alpha =0$ we obtain the approximation of zero forcing and dissipation, i.e. the``conservative form'' of the Lorenz-Krishnamurthy model studied by Vanneste \cite{24a}. Then, still using the \textit{Flow Curvature Method}, we will provide the thirteenth-order approximation of the \textit{slow manifold} associated with this model the invariance of which will be established according to Darboux theorem. Moreover, by posing $\varepsilon = 0$ in our \textit{slow manifold} analytical equation we will find again the first-order approximation of the \textit{slow manifold} given by Vanneste \cite{24a}.

\newpage

\subsection{The conservative LK model}

Thus, by using the same variables changing as previously and by posing $\delta = \varepsilon$ and $b = \epsilon$ in system (\ref{eq27}), Vanneste \cite{24a} obtained the following \textit{two-time scales singularly perturbed dynamical system}:

\begin{equation}
\label{eq34}
\left\{ {\begin{array}{*{20}c}
\varepsilon \dot{x} =  -y \\
\varepsilon \dot{y} = x + b u v\\
\dot{u} =  -v w + b \varepsilon v y\\
\dot{v} =  u w - b \varepsilon u y\\
\dot{w} = - u v
\end{array} } \right.
\end{equation}

where parameters $b$ and $\varepsilon$ control the strength of the coupling and the gravity-wave frequency, $x$ and $y$ are \textit{fast} modes while $u$, $v$ and $w$ are the \textit{slow} modes.

At the beginning of his paper Vanneste \cite{24a} gives the zero-order approximation in $\varepsilon$ (\textit{singular approximation}) of the slow manifold associated with the \textit{conservative LK model} (\ref{eq34}) by posing $\varepsilon = 0$:

\begin{equation}
\label{eq35}
\begin{aligned}
x & = - b u v,\\
y & = 0.
\end{aligned}
\end{equation}

Thus, as previously noticed by Lorenz \cite{15a} and recalled by Camassa \cite{2aa} and Vanneste \cite{24a} this model has an invariant manifold the equation of which is:

\begin{equation}
\label{eq36}
u^2 + v^2 = Cte
\end{equation}

Now, by using the \textit{Flow Curvature Method}, i.e. according to Prop. 1 the \textit{slow manifold} equation (\ref{eq9}) associated with \textit{conservative LK model} (\ref{eq34}) reads\footnote{See http://ginoux.univ-tln.fr for complete equation}:

\begin{equation}
\label{eq37}
\phi(\vec{X}, \varepsilon) = \det(\dot {\vec {X}},\ddot {\vec {X}},\dddot {\vec {X}},\ddddot {\vec {X}},\mathop{\vec {X}}\limits^{\ldots..} ) = 0
\end{equation}

As previously, it can be verified that the time derivative of the functional jacobian matrix of the \textit{fast} system (\ref{eq26}) (from which the \textit{slow} system (\ref{eq34}) has been deduced) is a zero matrix when $\varepsilon \to 0$. So, from \textit{Darboux Invariance Theorem} we can conclude that in the $\varepsilon$-vicinity of the \textit{singular approximation} the \textit{slow manifold} is invariant.

The implicit equation (\ref{eq37}) is a polynomial of degree $9$ for $u$, $v$ and $w$, of degree $5$ for $x$ and $11$ for $y$ and represents the thirteenth-order approximation in $\varepsilon$ of the \textit{slow manifold} of the \textit{conservative LK model} (\ref{eq34}).

By posing $\varepsilon = 0$ in the above Eq. (\ref{eq37}) we find that:

\begin{equation}
\label{eq38}
\phi(\vec{X}, 0) = (u^2 - w^2)(v^2 + w^2)((x + buv)^2 +y^2) = 0
\end{equation}

This slow manifold is made of a product of \textit{invariant manifolds} as it is easy to verify according to \textit{Darboux Invariance Theorem}. Let's compute the Lie derivative of the first and second term when $\varepsilon \to 0$ we have:

\begin{equation}
\label{eq39}
\begin{aligned}
L_{\vec X } ( u^2 - w^2 ) & = 0 \\
L_{\vec X } ( v^2 + w^2 ) & =  0\\
\end{aligned}
\end{equation}

Let's notice that the third term of Eq. (\ref{eq38}) is nothing else but the zero-order approximation in $\varepsilon$ (\textit{singular approximation}) of the \textit{slow manifold} (see Eq. (\ref{eq35})) given by Vanneste \cite{24a} which is also invariant when $\varepsilon \to 0$.

As previously, the decomposition into \textit{fast} and \textit{slow} modes enables to define a three-dimensional submanifold of the state space parametrized by $(u,v,w)$ and that he called \textit{slow manifold}. So, let's pose $x \to 0$ and $y \to 0$ in the above Eq. (\ref{eq37}) we find that:

\begin{equation}
\label{eq40}
\phi(u,v,w, \varepsilon) = (u^2 + v^2)(u^2 w^2 - u^4 + \varepsilon^2 v^2 w^2 (w^2 - (1 + b^2) u^2) ) = 0
\end{equation}

In addition to the \textit{quadratic invariant manifolds} (\ref{eq36}-\ref{eq39}) highlighted above we find another manifold. Let's compute its Lie derivative when $\varepsilon \to 0$ we obtain:

\begin{equation}
\label{eq41}
L_{\vec X } (u^2 w^2 - u^4 + \varepsilon^2 v^2 w^2 (w^2 - (1 + b^2) u^2) ) = 2 u v w ( u^2 - w^2 )
\end{equation}

Thus, we deduce that this manifold is \textit{locally} invariant, i.e. is invariant in the vicinity of the manifold defined by $u^2 - w^2  = 0$.

\newpage

The \textit{slow manifold} implicit equation (\ref{eq40}) associated with the \textit{conservative LK model} (\ref{eq34}) has been plotted in Fig. 1 in the ($u,v,w$) phase-space. Numerical integration of this model with a set of initial conditions ($x_0,y_0,u_0,v_0,w_0$) = ($2,2,-2,1.97,2$) taken on this \textit{slow manifold} (in blue on Fig. 1) enables to highlight that the \textit{trajectory curves} (in red on Fig. 1) ``visit'' every part of this hypersurface and stay in its $\varepsilon$-vicinity. The fixed point located at the origin has been plotted in green in the center of this figure.

\begin{figure}[htbp]
\centerline{\psfig{file=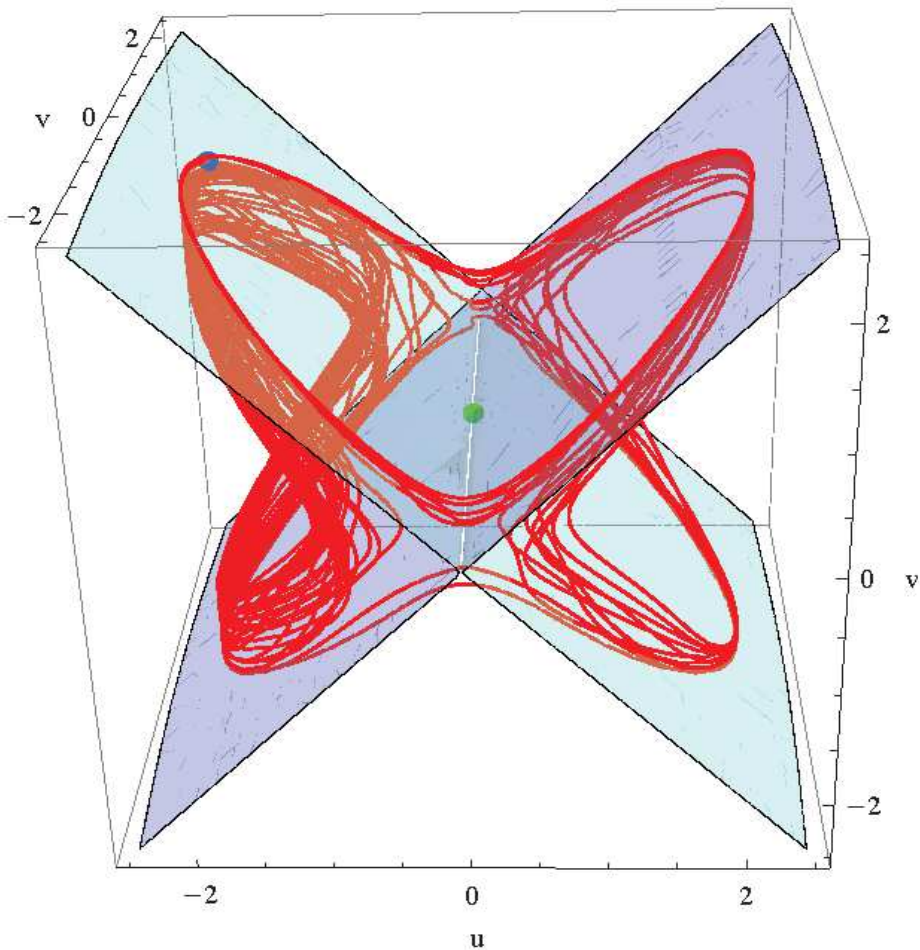,width=13cm,height=13cm}}
\caption{The \textit{conservative LK model slow invariant manifold} in $(u,v,w)$-space}
\label{fig3}
\end{figure}

\section{Discussion}

In this work the \textit{Flow Curvature Method} has enabled to provide the eighteenth-order approximation of the \textit{slow manifold} of the \textit{generalized LK model} and the thirteenth-order approximation of the \textit{conservative LK model} the invariance of which has been stated according to \textit{Darboux invariance theorem}.

\section*{Acknowledgements}

Authors would like to thank Prof. Jaume Llibre for his mathematical remarks and comments.

\newpage

\section*{APPENDIX}

The identity involved in the proof of the invariance of the \textit{slow manifold} (Sec. 4.2.2) is stated in this appendix.

\begin{align}
\label{eqA1}
& J\vec {a}_1 .\left( {\vec {a}_2 \wedge \vec {a}_3\wedge \ldots \wedge \vec {a}_n }
\right)+\vec {a}_1 .\left( {J\vec {a}_2 \wedge\vec {a}_3 \wedge \ldots \wedge \vec {a}_n }
\right)\hfill \notag \\
& +\ldots +\vec {a}_1 .\left( {\vec {a}_2 \wedge \vec {a}_3 \wedge \ldots \wedge J\vec
{a}_n } \right) = Tr\left( J \right)\vec {a}_1 .\left( {\vec {a}_2 \wedge
\ldots \wedge \vec {a}_n } \right)
\end{align}

\textit{\textbf{Proof.}} The proof is based on inner product properties.

To the functional jacobian matrix $J$ is associated an eigenbasis: $\left\{ \vec {Y_{\lambda _1 }} ,\vec {Y_{\lambda _2 }}, \ldots,\vec {Y_{\lambda _n }} \right\}$.

Let suppose that there exists a transformation\footnote{By considering that each vector $\vec {a}_i $ may be spanned on the eigenbasis, calculus is longer and tedious but leads to the same result.} such that:\\
to each vector $\vec {a}_i $ corresponds the eigenvector $\vec {Y_{\lambda _i }}$
with $i=1,\ldots,n$.

Each inner product of the left hand side Eq. (\ref{eqA1}) may be transformed into

\begin{align*}
& J\vec {a}_1 \cdot \left( {\vec {a}_2 \wedge \vec {a}_3 \wedge \ldots \wedge \vec
{a}_n } \right)=\lambda _1 \vec {a}_1 \cdot \left( {\vec {a}_2 \wedge \vec
{a}_3 \wedge \ldots \wedge \vec {a}_n } \right)=\lambda _1 \vec {a}_1 \cdot
\left( {\vec {a}_2 \wedge \vec {a}_3 \wedge \ldots \wedge \vec {a}_n } \right)\\
& \vec {a}_1 \cdot \left( {J\vec {a}_2 \wedge \vec {a}_3 \wedge \ldots \wedge \vec
{a}_n } \right) = \vec {a}_1 \cdot \left( {\lambda _2 \vec {a}_2 \wedge \vec
{a}_3 \wedge \ldots \wedge \vec {a}_n } \right) = \lambda _2 \vec {a}_1 \cdot
\left( {\vec {a}_2 \wedge \vec {a}_3 \wedge \ldots \wedge \vec {a}_n } \right)\\
& \hfill {\ldots}{\ldots}{\ldots}{\ldots}{\ldots}{\ldots}{\ldots}{\ldots}{\ldots}{\ldots}{\ldots}{\ldots}{\ldots}{\ldots}{\ldots}{\ldots}{\ldots}{\ldots}{\ldots}{\ldots}{\ldots}{\ldots}{\ldots}{\ldots}{\ldots}{\ldots}{\ldots}{\ldots}{\ldots}{\ldots}\hfill \\
& \vec {a}_1 \cdot \left( {\vec {a}_2 \wedge \vec {a}_3 \wedge \ldots \wedge J\vec
{a}_n } \right) = \vec {a}_1 \cdot \left( {\vec {a}_2 \wedge \vec {a}_3 \wedge
\ldots \wedge \lambda _n \vec {a}_n } \right) = \lambda _n \vec {a}_1 \cdot \left(
{\vec {a}_2 \wedge \vec {a}_3 \wedge \ldots \wedge \vec {a}_n } \right)
\end{align*}

Making the sum of these factors the proof is stated.


\begin{thebibliography}{00}


% notes:

% \bibitem{label2}
% If there is a note, it should come last:
% \bibitem{label3} \note
% \end{subbibitems}

\bibitem{1} A.A. Andronov {\&} S.E. Chaikin: Theory of Oscillators, I. Moscow (1937); English transl., Princeton Univ. Press, Princeton, N.J. (1949)

\bibitem{1a} J.P. Boyd: The slow manifold of five-mode model. J. Atmos. Sci. 51 (8) 1057--1064 (1994)

\bibitem{2aa} R. Camassa: On the geometry of an atmospheric slow manifold. Physica D 84, 357--397 (1995)

\bibitem{2a} R. Camassa \& Siu-Kei-Tin: The global geometry of the slow manifold in the Lorenz–
Krishnamurthy model. J. Atmos. Sci. 53, 3251--3264 (1996)

\bibitem{2} J.D. Cole, Perturbation Methods in Applied Mathematics, Blaisdell, Waltham, MA (1968)

\bibitem{3} G. Darboux: Mémoire sur les équations différentielles algébriques du premier ordre et du premier degré. Bull. Sci. Math. S\'{e}r. 2(2), 60--96, 123--143, 151--200 (1878)

\bibitem{5} N. Fenichel: Persistence and Smoothness of Invariant Manifolds for Flows. Ind. Univ. Math. J. 21, 193--225 (1971)

\bibitem{6} N. Fenichel: Asymptotic stability with rate conditions. Ind. Univ. Math. J. 23, 1109--1137 (1974)

\bibitem{7} N. Fenichel: Asymptotic stability with rate conditions II. Ind. Univ. Math. J. 26, 81--93 (1977)

\bibitem{8} N. Fenichel: Geometric singular perturbation theory for ordinary differential equations, J. Diff. Eq. 31, 53--98 (1979)

\bibitem{8a} A.C. Fowler {\&} G. Kember: The Lorenz–Krishnamurthy slow manifold. J. Atmos. Sci. 53 (10) 1433--1437  (1996)

\bibitem{9} J.M. Ginoux {\&} B. Rossetto: Slow manifold of a neuronal bursting model. In Emergent Properties in Natural and Articial Dynamical Systems (eds) M.A. Aziz-Alaoui \& C. Bertelle, pp. 119--128. Springer-Verlag, Heidelberg (2006)

\bibitem{9a} J.M. Ginoux {\&} B. Rossetto: Differential Geometry and Mechanics Applications to Chaotic Dynamical Systems. Int. J. Bif. {\&} Chaos 4(16) 887--910 (2006)

\bibitem{9b} J.M. Ginoux, B. Rossetto \& L.O. Chua: Slow Invariant Manifolds as Curvature of the Flow of Dynamical Systems. Int. J. Bif. {\&} Chaos 11(18) 3409--3430 (2008)

\bibitem{9c} J.M. Ginoux, Differential geometry applied to dynamical systems, World Scientific Series on Nonlinear Science, Series A {\bf 66} World Scientific, Singapore (2009)

\bibitem{9d} J.M. Ginoux \& J. Llibre: The flow curvature method applied to canard explosion. J. Phys. A: Math. Theor. 44 (2011). 465203 doi:10.1088/1751-8113/44/46/465203

\bibitem{9e} J.M. Ginoux, J. Llibre \& L.O. Chua: Canards from Chua's circuit. Int. J. Bif. {\&} Chaos 23(4) (2013). 1330010 doi:10.1142/S0218127413300103

\bibitem{10b} J. Guckhenheimer {\&} Ph. Holmes, Nonlinear Oscillations, Dynamical Systems and Bifurcation of Vector Fields, Springer-Verlag, New York, (1983)

\bibitem{12} M. W. Hirsch, C.C. Pugh {\&} M. Shub: Invariant Manifolds, Springer-Verlag, New York (1977)

\bibitem{12a} S.J. Jacobs: Existence of a slow manifold in a model system of equations. J. Atmos. Sci. 48, 893--901 (1991)

\bibitem{13} C.K.R.T. Jones: Geometric Singular Pertubation Theory. In Dynamical Systems, Montecatini Terme, L. Arnold, Lecture Notes in Mathematics, vol. 1609, pp 44--118, Springer-Verlag, Heidelberg (1994)

\bibitem{14} T. Kaper: An Introduction to Geometric Methods and Dynamical Systems Theory for Singular Perturbation Problems. In Analyzing
multiscale phenomena using singular perturbation methods, (Baltimore, MD, 1998), pp 85--131, Amer. Math. Soc., Providence, RI (1999)

\bibitem{15} N. Levinson: A second-order differential equation with singular solutions. Ann. Math. 50, 127--153 (1949)

\bibitem{15aaa} C. E. Leith: Nonlinear Normal Mode Initialization and Quasi-Geostrophic Theory. J. Atmos. Sci. 37, 958--968 (1980)

\bibitem{llibre} J. Llibre, R. Saghin {\&} X. Zhang: On the analytic integrability of the 5-dimensional Lorenz system for the
gravity-wave activity. To appear in Proc. Amer. Math. Soc.

\bibitem{15aa} E.N. Lorenz: Attractor sets and quasi-geostrophic equilibrium. J. Atmos. Sci. 37,  1685--1699 (1980)

\bibitem{15a} E.N. Lorenz: On the existence of a slow manifold. J. Atmos. Sci. 43, 1547--1557 (1986)

\bibitem{15b} E.N. Lorenz {\&} V. Krishnamurthy: On the nonexistence of a slow manifold. J. Atmos. Sci. 44, 2940--2950 (1987)

\bibitem{15c} E.N. Lorenz: The slow manifold — What is it? J. Atmos. Sci. 49, 2449--2451 (1992)

\bibitem{16} R.E. O'Malley, Introduction to Singular Perturbations, Academic Press, New York (1974)

\bibitem{17} R.E. O'Malley, Singular Perturbations Methods for Ordinary Differential Equations, Springer-Verlag, New York (1991)

\bibitem{17a} M. Phani Sudheer, Ravi S. Nanjundiah {\&} A.S. Vasudeva Murthy: Revisiting the slow manifold of the Lorenz-Krishnamurthy quintet, Discrete and Continuous Dynamical Systems Serie B, Vol. 6, N\textsuperscript{o} 6, 1403--1416 (2006)

\bibitem{21} B. Rossetto, T. Lenzini, S. Ramdani {\&} G. Suchey: Slow-fast autonomous dynamical systems. Int. J. Bif. {\&} Chaos 8(11) 2135--2145 (1998)

\bibitem{23} A.N. Tikhonov: On the dependence of solutions of differential equations on a small parameter. Mat. Sbornik N.S., 31, 575--586 (1948)

\bibitem{24} B. Van der Pol: On 'Relaxation-Oscillations'. Phil. Mag., 7(2) 978--992 (1926)

\bibitem{24a} J. Vanneste: Asymptotics of a slow manifold. SIAM J. Appl. Dynam. Syst. 7, 1163--1190 (2008)

\bibitem{25} W.R. Wasow: Asymptotic Expansions for Ordinary Differential Equations, Wiley-Interscience, New York (1965)


\end{thebibliography}
\end{document}